\documentclass[12pt]{article}
 \usepackage{amsfonts}
 \usepackage{amssymb}
 \usepackage{amsmath,amsxtra}

\newtheorem{theorem}{Theorem}

\newtheorem{proposition}{Proposition}

\newtheorem{lemma}{Lemma}

\newtheorem{cor}{Corollary}

\newtheorem{rem}{Remark}

 \numberwithin{equation}{section}

\newcommand{\npa}{\addtocounter{num}{1} \noindent
{\bf \arabic{section}.\arabic{num}}\;\;}
\newcommand{\proof}{{\bf Proof.}~}

\newcommand{\bqa}{\begin{eqnarray}}
\newcommand\eqa {\end{eqnarray}}
\newcommand{\beq}{\begin{eqnarray}}
\newcommand{\beqn}{\begin{eqnarray}\nonumber}
\newcommand{\eeq}{\end{eqnarray}}
\newcommand{\be}{\begin{array}}
\newcommand{\ee}{\end{array}}

 \newcommand{\pr}{\partial}

 \newcommand{\rk}{\mathrm{rk}}

 \newcommand{\Ker}{\mathrm{Ker}}

 \newcommand{\cA}{{\cal A}}
 
 \newcommand{\cM}{{\cal M}}
 
 \newcommand{\cO}{{\cal O}}
 \newcommand{\cB}{{\cal B}}

 \newcommand{\cX}{{\cal X}}

 \newcommand{\cN}{{\cal N}}

 \newcommand{\R}{{\mathbb R}}
 \newcommand{\Z}{{\mathbb Z}}

   \newcommand{\g}{{\mathfrak g}}
   \newcommand{\gl}{{\mathfrak gl}}

   \newcommand{\sll}{{\mathfrak sl}}

   \newcommand{\h}{{\mathfrak h}}

\begin{document}

   \def\sp{\mathfrak sp}
   \def\sll{\mathfrak sl}
   \def\P{{\mathbb P}}
   \def\H{\mathbb H}
   \def\a{\alpha}
   \def\b{\beta}
   \def\t{\theta}
   \def\la{\lambda}
   \def\e{\epsilon}

 \def\gr{\g^{\scriptscriptstyle\mathrm{gr}}}
 \def\godd{\g_{\scriptscriptstyle 1}}
 \def\geven{\g_{\scriptscriptstyle 0}}
 \def\grodd{\gr_{\scriptscriptstyle 1}}
 \def\greven{\gr_{\scriptscriptstyle 0}}

  \def\sst{\scriptscriptstyle}

   %\today
\vskip 10mm

\title{%The Lie algebra of a supermanifold
Lie algebaic characterization of\\ supercommutative space}

\author{Janusz Grabowski\thanks{The
research of J. Grabowski was financed by the Polish
Ministry of Science and Higher Education %, by means of the budget for science 2006-2009,
under the grant No. N201 005 31/0115.} , Alexei Kotov, Norbert
Poncin\thanks{The research of N. Poncin was supported by UL-grant
SGQnp2008.}}

\date{}

\maketitle

\abstract{During the last decades algebraization of space turned
out to be a promising tool at the interface between Mathematics
and Theoretical Physics. Starting with works by
Gel'fand-Kolmogoroff and Gel'fand-Naimark, this branch developed
as from the fortieth in two directions: algebraic characterization
of usual geometric space on the one hand, and algebraically
defined noncommutative space, which is known to be tightly related
with e.g. quantum gravity and super string theory, on the other
hand. In this note, we combine both aspects, prove a superversion
of Shanks and Pursell's classical result stating that any
isomorphism of the Lie algebras of compactly supported vector
fields is implemented by a diffeomorphism of underlying manifolds.
We thus provide a super Lie algebraic characterization of super
and graded spaces and describe explicitly isomorphisms of the
super Lie algebras of super vector fields.} \vskip 4mm

 %\vspace{0.2em}
 \vskip 2mm


\section{Introduction}

{\it Algebraic characterization of space} can be traced back to
Gel'fand and Kolmogoroff, who proved in 1939 that two compact
topological spaces are homeomorphic if and only if the algebras of
continuous functions growing on them are isomorphic. A similar
result for second countable smooth manifolds and the algebras of
smooth functions is also regarded as classical (for the general
case see \cite{Gra05, Mrc05}). In 1954, Pursell and Shanks
\cite{Pursell-Shanks} substituted the Lie algebra of compactly
supported vector fields of a manifold for the commutative
associative algebra of smooth functions. This classical upshot
triggered a multitude of papers on similar issues by many
different authors, which we extensively depicted in our previous
works. In 2004, two of us proved Pursell-Shanks type results for
the Lie algebra of differential operators of a manifold, and for
the Poisson-Lie algebra of smooth functions on the cotangent
bundle that are polynomial along the fibers. Our results indicate
once more a ``no-go'' theorem for the Dirac quantization problem,
as they imply that the preceding Lie algebras are not isomorphic
-- since they have nonisomorphic automorphism groups. Let us
mention, hoping that the remark might instigate further progress,
that the last observation is tightly related to the
Kanel-Kontsevich conjecture that maintains that the automorphism
groups of the Weyl algebra -- modelled on the algebra of
differential operators {\it with polynomial coefficients} -- and
of the standard {\it Poisson algebra of polynomials} are
isomorphic!\medskip

Another landmark in the field of {\it algebraization of space} is
the Gel'fand-Naimark theorem, 1943, which states that any
$C^*$-algebra is isometrically $*$-isomorphic to a $C^*$-algebra
of bounded operators on a Hilbert space. This result is usually
viewed as the starting point of noncommutative geometry: the basic
idea of this branch is to treat certain noncommutative algebraic
structures that arise in Physics as if they were related to some
``noncommutative spaces'', although there are no such spaces in
the usual sense of the word. It is well-known that {\it
algebraically defined noncommutative space} is an important tool
in quantization of gravity, i.e. in the attempt to unify the
contradictory concepts of gravity (which makes no sense at
0-distance) and quantum theory (which precisely concerns the
``0-distance''). Indeed, one of the ways out of this conflicting
situations consists in replacing at small distance the usual
commutative space by noncommutative space. Another possible
remedy, which allows dealing with the mentioned singularities, is
the replacement of points by extended geometric objects or
strings, viewed as the fundamental constituents of reality. The
effort to incorporate fermions, the building blocks of matter, in
the spectrum of string theory led to supersymmetry and superspace
(resp. $\Z$-graded space) -- a particular type of noncommutative
space: a supermanifold (resp. $\Z$-graded manifold) is a sheaf of
supercommutative (resp. $\Z$-graded commutative) associative
algebras that is locally isomorphic with a free Grassmann algebra
with coefficients in the functions of Euclidean space.\medskip

In the present note, we combine the two aforementioned aspects of
algebraization of space --  algebraic characterization of usual
space and algebraically defined noncommutative space. More
precisely, we prove that if two Lie superalgebras of supervector
fields are isomorphic if and only if the underlying smooth
supermanifolds are diffeomorphic as sheafs of supercommutative
algebras and describe explicitly automorphism groups of the super
Lie algebras of super vector fields.\medskip

The paper is organized as follows. In section 2, we recall that
every (smooth) supermanifold ${\cal M}$ is (noncanonically)
diffeomorphic to the total space of some vector bundle $V$ with
reversed parity in the fibres \cite{Gawedski}, which we denote by
$\Pi V$ and which is actually the prototype of a $\Z$-graded
manifold. Further, we show that the super Lie algebra of vector
fields of ${\cal M}$ admits a canonical Lie-algebraic filtration,
such that the corresponding quotient is isomorphic to the
$\Z$-graded Lie algebra of vector fields of $\Pi V$, whatever
diffeomorphism is chosen. In section 3, we prove the
aforementioned superversion of Shanks and Pursell's classical
result in the case of $\Z$-graded manifolds $\Pi V$, and finally,
in section 4, we use the results of section 2 to deduce the
supercase from the preceding $\Z$-graded case.

\section{Lie algebra of supervector fields}

\npa Let $\cM$ be a smooth supermanifold of dimension $(s,r)$ over the body $M$. Here we understand the
supermanifold as a ringed space: the standard manifold $M$ of dimension $s$ is equipped with a sheaf $\cO_M$
of superalgebras which is locally isomorphic to $C^\infty(\R^{s})\otimes\Lambda^\bullet(\xi^1,\dots,\xi^{r})$.
Sections of this sheaf form a superalgebra $\cA=\cA_{\sst 0}\oplus\cA_{\sst 1}$ of {\it smooth functions} on
the supermanifold $\cM$.

An important result of smooth supergeometry \cite{Gawedski} (see also \cite{Manin_Gauge}) asserts that there
exists a vector bundle $V$ of rank $r$ over $M$, such that $\cM$ is diffeomorphic as a supermanifold to $\Pi
V$, that is, to the total space of $V$ with the reversed parity of fibres. This implies that the algebra of
smooth functions on $\cM$ is isomorphic (as a commutative superalgebra) to the algebra of functions on $\Pi
V$, which is canonically identified with $\Gamma (\Lambda^\bullet V^*)$. This isomorphism is not canonical but
it gives us an identification
 \beq\label{isomorphism_functions}
  \cA=\cA_{\sst 0}\oplus\cA_{\sst 1}\simeq \Gamma (\Lambda^\bullet V^*)\,,
 \eeq
 with
 \beq\label{isomorphism_functions1}
  \cA_{\sst 0}=\bigoplus\limits_{i\ge 0} \cA^{\sst 2i}\,,\hspace{2mm}
  \cA_{\sst 1}=\bigoplus\limits_{i\ge 0} \cA^{\sst 2i+1}\,,\hspace{2mm}
  \mathrm{where}\hspace{1mm}\cA^{\sst k}\simeq \Gamma (\Lambda^{\sst k} V^*)\,.
 \eeq

\vskip 2mm\noindent The choice of an isomorphism
(\ref{isomorphism_functions}) provides therefore an additional $\Z-$grading
in $\cA$ which is compatible with the given super-structure. Such a grading
uniquely determines the {\it Euler vector
field}, that is, an operator $\e$ satisfying the following property (the
definition of $\e$ implements the Leibnitz rule):
 \beq\label{Euler_field}
  \cA^{\sst k}=\{a\in\cA \,\mid \, \e (a)=k a\}\,.
 \eeq
Denote with $\g$ the super Lie algebra of vector fields on $\cM$, the even and odd parts of which are $\geven$
and $\godd$, respectively,
 \beqn
   \g=\geven\oplus\godd\;.
 \eeq
 To the end of this section an isomorphism \ref{isomorphism_functions}
is chosen. Moreover, we will assume that the rank $r$ of the vector bundle $V$ is at least 1, otherwise we are
in the standard purely even situation.
\begin{proposition} The adjoint action of the Euler vector field
supplies $\g$ with a $\Z$-grading compatible with the Lie super
structure such that
 \beq
 \geven=\bigoplus\limits_{i\ge 0} \g^{\sst 2i}\,,\hspace{2mm}
 \godd=\bigoplus\limits_{i\ge -1} \g^{\sst 2i+1}\,,\hspace{2mm}
  \mathrm{where}\hspace{1mm}\g^{\sst k}\colon=\{ X\in \g\, \mid\, [\e, X]=kX\}\,.
 \eeq
Any super vector field $X$ admits a unique homogeneous
decomposition $X=\sum\limits_{m\ge -1}X_m$ with respect to the
Euler vector field, $[\e, X_m]=mX_m$. In local coordinates this
decomposition is given by the polynomial degree of $\xi^a$:
 \beq\label{vf_local_decomposition}
X_m:=\sum\limits_{a_1<\ldots < a_m}\sum\limits_i f_{a_1\ldots
a_m}^i (x)\xi^{a_1}\ldots \xi^{a_m} \frac{\pr}{\pr x^i}
+\\\nonumber +\sum\limits_{b_1<\ldots < b_{m+1}}\sum\limits_c
  g_{b_1\ldots  b_{m+1}}^c (x)\xi^{b_1}\ldots \xi^{b_{m+1}}
\frac{\pr}{\pr \xi^c}
  \;.
 \eeq
Here $f_{b_1\ldots  b_r}^i (x) $ and $g_{a_1\ldots  a_r}^b (x)$
are smooth functions of $x$.
 \end{proposition}

\noindent\proof Suppose we are given a local trivialization of $V$ over $M$, that is, an open cover by
coordinate charts $\{U_\a, x^i\}$ together with a local frame of the restriction of $V$ to each $U_\a$,
denoted by $(e_a)$, where $a=1,\ldots, r$. Combining these data, we obtain a local coordinate description of
$\cM$: $\{U_\a, x^i, \xi^a\}$, where $\xi^a$ are dual to $e_a$ thought of as odd coordinates. By construction
of the trivialization, the change of coordinates over double overlaps is linear with respect to odd
coordinates, so the aforementioned Euler vector field is well defined, and in any coordinate system reads
 \beq\label{Euler_vf_formula}
 \e =\sum\limits_a \xi^a \frac{\pr}{\pr \xi^a}\,.
 \eeq
It is now a standard task in local coordinates to write any vector
field in the form (\ref{vf_local_decomposition}). $\square$

\begin{rem}{\rm
Apparently, each vector field of degree -1 can be identified with
a section of $V$, thus $\g^{-1}$ is naturally isomorphic to
$\Gamma (V)$ which acts on
 $\cA$ by contractions, provided we identify $\cA$ with
 $\Gamma (\Lambda^\bullet V^*)$. On the other hand, a super vector field of degree 0, which can be written as
 \beq\label{degree_zero_vf}
  X=\sum\limits_i f^i (x)\frac{\pr}{\pr x^i} +
  \sum\limits_{a,b}
  g_{a}^b (x)\xi^{a}
\frac{\pr}{\pr \xi^b}\;,
 \eeq
defines a general infinitesimal automorphism of the vector bundle
$V$. The vector fields from $\g^{\sst 0}$ can be therefore
identified with the sections of the Lie algebroid of infinitesimal
automorphism of $V$, called sometimes the {\it  Atiyah algebroid}
of $V$. This identification respects the bracket, i.e. is a Lie
algebra isomorphism. The corresponding anchor map $\rho:\g^{\sst
0}\rightarrow\cX(M)$ from the Lie algebra $\g^{\sst 0}$ of the
Atiyah algebroid into the Lie algebra $\cX(M)$ of vector fields on
$M$ in local coordinates reads
 \beq\label{anchor}
  \rho(X)=\rho\left(\sum\limits_i f^i (x)\frac{\pr}{\pr x^i} +
  \sum\limits_{a,b}  g_{a}^b (x)\xi^{a}\frac{\pr}{\pr \xi^b}\right)
=\sum\limits_i f^i (x)\frac{\pr}{\pr x^i}
 \eeq
 and is also a Lie algebra homomorphism.}
\end{rem}

\vskip 2mm\npa Suppose we are given a maximal ideal of $\geven$,
denoted by $\g'$, the elements of which act as $ad$-nilpotent
operators in $\g$, that is, for each $X\in\g'$ there exists a
non-negative integer $m$, such that $ad_X^m (Y)=0$ for all $Y\in
\g$.

 \begin{proposition}\label{maximal_nilpotent}
 The ideal $\g'$ is related to the $\Z-$grading by the
 formula:
  \beq
   \g' =\bigoplus\limits_{i>0}\g^{\sst 2i}\,.
  \eeq
  Since the ideal $\g'$ is defined in purely super Lie algebraic terms, the
  latter space in fact does no depend on the introduced $\Z-$graduation.
 \end{proposition}

\noindent\proof Let us consider the image of $\g'$ under the
 projection
\beqn
 \pi_{\sst 0}\colon \geven\to \geven /\oplus_{i>0}\g^{\sst 2i}\simeq
\g^{\sst 0}\,,
 \eeq
 which is fixed by the choice of an isomorphism $\cM\simeq \Pi V$
 (\ref{isomorphism_functions}).
 Apparently, $\pi_{\sst 0}(\g')$ is a maximal
$ad-$nilpotent ideal of $\g^{\sst 0}$. Let us apply the anchor map
$\rho$ to this ideal. For each $X\in \pi_{\sst 0}(\g')$, the
vector field $\rho (X)$ has to be $ad-$nilpotent as well. It is
easy to see that such a vector field is necessarily zero, as any
vector field on a standard (even) manifold can be written locally
as the coordinate vector field $\partial_{x_1}$ in a neighborhood
of any point at which it does not vanish.

\vskip 2mm\noindent Hence we conclude that $\pi_{\sst
0}(\g')\subset \Ker \rho$. But $\Ker \rho$ is the bundle of Lie
algebras over $M$, the fiber of which is isomorphic to $\gl (V_z)$
at any $z\in M$. Thus $\pi_{\sst 0}(\g')$ evaluated at $z$ is
necessarily a subset of the ideal of scalar operators,  the only
ideal of  $\gl (V_z)$, which implies that $\pi_{\sst
0}(\g')\subset A\e$. Here $\e$ is the Euler vector field and
$A\colon =\cA^{\sst 0}$ is the algebra of functions on $M$.

\vskip 2mm\noindent Let us assume that $X\in\g_0$, then $\pi_{\sst
0} (X)=f\e$ for some smooth function $f\in A$. Therefore $X=f\e+$
terms of order $\ge 2$. On the other hand, $[f\e, Y]=-fY$ for each
$f\in A $ and $Y\in \g^{\sst -1}$, hence $ad_X^m (Y)=(-f)^m Y+$
terms of order $\ge 1$. Thus we conclude that $X$ is
$ad-$nilpotent if and only if $f=0$, therefore $\pi_{\sst 0}(X)=0$
and $\g'\subset \Ker\pi_{\sst 0}$. The kernel of $\pi_{\sst 0}$ is
an ideal of $\geven$, consisting of $ad-$nilpotent elements; but
$\g'$ has to be maximal, which immediately implements the identity
$\g'=\Ker\pi_{\sst 0}$. $\square$\medskip

\vskip 2mm\noindent There is an additional conclusion drawn from the above proof which
we formulate as a separate proposition that we will use later.
\begin{proposition}\label{pf}
The maximal Lie ideal in $\g^0$ of elements acting
$ad$-nilpotently on $\g^0$ consists of vector fields of the form
$f\epsilon$ with $f\in A$ being a smooth function on the body $M$.
\end{proposition}

\vskip 2mm\noindent Let us introduce (inductively) the following
subspaces:
 \beq
  \g^{\sst (p+2)}\colon =[\g', \g^{\sst (p)}]\,,\hspace{3mm}
  \mathrm{where}\hspace{1mm} \g^{\sst (-1)}\colon =\godd \,,
  \hspace{1mm} \g^{\sst (0)}\colon =\geven \,.
 \eeq

\begin{proposition}\label{maximal_nilpotent2}
Let us assume that $r=\rk V> 2$ or $\dim M>0$ and $\rk V>1$
  Then
  \beq
   \g^{\sst (p)}=\bigoplus\limits_{i\ge 0} \g^{\sst p+2i}\,.
  \eeq
 independently on the choice of the isomorphism
 (\ref{isomorphism_functions}).
 \end{proposition}
 \noindent\proof It
follows from the equality $[\g^{\sst p}, \g^{\sst q}]= \g^{\sst
p+q}$, which is true for all $p,q$ (except for $p=q=0$ in the pure
odd case when $M$ is a point). Indeed, let $X$ be a super vector
field of degree $p+q$, then given a local cover $\{U_\a\}$ of $M$,
we decompose $X$ into a sum of $X_\a$, such that $\mathrm{supp}
X_\a\subset M$. It is enough to find $Y_\a^k$, $Z_\a^k$ of the
degree $p$ and $q$, respectively, with the support in $U_\a$ for
each $\a$, such that $X_\a=\sum_k [Y_\a^k, Z_\a^k]$. Now we use
the local representation (\ref{vf_local_decomposition}) of $X$.
$\square$\medskip

\vskip 2mm\noindent We still make the assumption that $r=\rk V> 2$
or $\dim M>0$ and $\rk V>1$, leaving the low rank cases, for which
the next two corollaries fail, to separate considerations.
\begin{cor} The filtration of $\g$ by $\g^{\sst (p)}$ respects the
super Lie algebra structure, i.e. $[\g^{\sst (p)}, \g^{\sst (q)}]\subset \g^{\sst (p+q)}$. The graded Lie
superalgebra $\gr$, associated with the given filtration,
 \beqn
\gr=\bigoplus\limits_{p\ge -1} (\gr)^{\sst p}\,, \hspace{3mm} \mathrm{where}\hspace{1mm} (\gr)^{\sst p} =
\g^{\sst (p)}/\g^{\sst(p+2)}\,,
 \eeq
equipped with the bracket naturally induced by the bracket in $\g$,
 is isomorphic to $\g$ as a $\Z-$graded superalgebra, independently
 on the choice of (\ref{isomorphism_functions}).
\end{cor}
Since the filtration $\g^{\sst (p)}$ is canonical, thus preserved by any automorphism of $\g$, any
automorphism $\psi$ of the Lie superalgebra $\g$ induces an automorphism $ \mathrm{p}(\psi)$ of the graded Lie
algebra $\gr$ by
\beq\label{mor}\mathrm{p}(\psi)([X])=[\psi(X)]\,,
\eeq
where $[X]$ is the coset of $X\in\g^{\sst (p)}$.

\begin{cor}\label{morphism_groups}
 The formula (\ref{mor}) defines a group homomorphism
 \beq\label{morphism_p}
 \mathrm{p}\colon Aut_{\Z_2} (\g)\to Aut_{\Z} (\gr)\,,
 \eeq
 where the former and the latter groups
 consist of all automorphisms preserving $\Z_2$ on $\g$ and
 $\Z-$grading on $\gr$, respectively.
\end{cor}

\section{Lie algebras of $\Z$-graded manifolds associated
with vector bundles}

Let $\cM=\Pi V$ and $\cN=\Pi W$ be $\Z$-graded manifolds
associated with vector bundles $V\to M$ and $W\to N$ respectively.
Let $\cA=\oplus_{i\ge 0}\Gamma(\Lambda^iV^\ast)$ and
$\cB=\oplus_{i\ge 0}\Gamma(\Lambda^iW^\ast)$ be the graded
algebras of smooth functions on $\cM$ and $\cN$ (Grassmann
algebras of multi-sections of dual bundles), respectively. Let
$\g$ and $\h$ be the corresponding $\Z$-graded Lie algebras of
vector fields. Let us also assume that $\dim M$ and $\dim N$ are
non-zero or the ranks $\rk V$ and $\rk W$ are both positive and
different from $2$.

\begin{theorem}\label{Isomorphism_Z}
For any isomorphism $\psi\colon \g \to\h$ of the $\Z$-graded Lie
algebras of vector fields on $\Pi V$ and $\Pi W$, respectively,
there exists an isomorphism of vector bundles $\phi \colon V\to W$
such that $\psi (X)=(\phi^*)^{-1}\circ X \circ\phi^*$, where
$\phi^*\colon\cB\to\cA$ is the isomorphism of the $\Z$-graded
algebras of smooth functions (multi-sections of the dual bundles)
induced by $\phi$.
\end{theorem}

\noindent\proof Let us restrict the isomorphism $\psi$ to
$\g^{\sst 0}$. According to proposition \ref{pf}, the subspaces
$A\e$ and $B\e$, where $A$ and $B$ are the algebras of functions
on $M$ and $N$, respectively, are the (uniquely determined)
maximal ideals of $\g^{\sst 0}$ and $\h^{\sst 0}$ acting
nilpotently on $\g^{\sst 0}$ and $\h^{\sst 0}$ respectively,
therefore $B\e =\psi (A\e)$. Of course, here we identify the two
Euler vector fields in our bundles as uniquely determined by the
$\Z-$graduation. This implies the existence of a bijective map
$\widetilde\psi \colon B\to A$, such that
 \beq
  \psi (f\e)=\widetilde\psi^{-1} (f)\e\,, \hspace{2mm} \forall f\in
  A\,.
 \eeq
Taking into account that $[X, f\e]=\rho (X)(f)\e$, where $\rho
(X)$ is the anchor of $X\in \g^{\sst 0}$, we immediately obtain
the following property
 \beq \rho (\psi (X))=\widetilde\psi^{-1} \rho
(X)\widetilde\psi\,,
 \eeq
which implies that the conjugation by $\widetilde\psi$ induces a
Lie algebra isomorphism $\rho (\g^{\sst 0})\to \rho (\h^{\sst
0})$. On the other hand, these Lie algebras consist of all vector
fields on $M$ and $N$, correspondingly. Using the classical result
on the Lie algebras of vector fields \cite{Grabowski} (see also
\cite{Pursell-Shanks,Amemiya}), we conclude that the conjugation
by $\widetilde\psi$ coincides with the conjugation by some
diffeomorphism $h: M\to N$. One can also conclude the latter fact
from a theorem in \cite{Grabowski-Grabowska} applied directly to
the Atiyah algebroid $\g^0$. Any operator, acting on smooth
functions and commuting with all vector fields, is necessarily a
constant, which implies that $\widetilde\psi (h^*)^{-1}$ is the
operator of multiplication by a non-zero constant $\lambda$, and
thus $\widetilde\psi =\la h^*$. But the Euler vector fields are
uniquely defined, so they are associated by the isomorphism which
yields $\lambda=1$.

\vskip 2mm\noindent Now we restrict the automorphism $\psi$ to
$\g^{\sst -1}$; thus we obtain a non-degenerate linear map $\Gamma
(V)\to \Gamma (W)$. If we proved that for each $s\in \Gamma (V)$
and $f\in A$, the following property holds: $\psi (fs)=(h^*)^{-1}
(f)\psi (s)$, the restriction of $\psi$ would have been induced by
a bundle map $V\to W$ covering $h$. Indeed, $fs=[s, f\e]$. On the
other hand, $\psi (f\e)=\widetilde\psi^{-1} (f)\e$, so that
 \beqn \psi (fs)=\psi [s, f\e]=[\psi (s), \psi
(f\e)]=[\psi (s), \widetilde\psi^{-1} (f)\e]=\widetilde\psi^{-1}
(f)\psi (s)\,.
 \eeq
 Since $\widetilde\psi$ is an automorphism of the algebra of
 functions, $\psi (fs)=(h^*)^{-1} (f)\psi (s)$. The bundle map $V\to W$,
induced by the restriction of $\psi$ to $\g^{\sst -1}$, can be
uniquely extended to a diffeomorphism $\phi$ of $\cM$ over $\Z$.
The uniqueness follows from the isomorphism
(\ref{isomorphism_functions}) because the algebras of functions on
$\cM$ and $\cN$ are freely generated by $\Gamma (V^*)$ and $\Gamma
(W^*)$, respectively. Taking into account that the representation
of $\g^{\sst 0}$ on the space of sections of $V$ is faithful, we
immediately obtain the required property of Theorem
\ref{Isomorphism_Z} for $\g^0\oplus \g^{-1}$.

\vskip 2mm \noindent The last step is proving that the extension
of the restriction of $\psi$ to $\g^{\sst 0}\oplus \g^{\sst -1}$
is unique. Suppose there exists another  Lie algebra morphism
$\psi'$, satisfying the property \beq\label{another_extension}
\psi (X)=\psi'(X)\, ,\hspace{3mm} \forall X\in \g^{\sst 0}\oplus
\g^{\sst -1}\,. \eeq Then for each $Z\in \g^{\sst k}$, where $k\ge
0$, and $X_i\in \g^{\sst -1}$, $i=1,\ldots, k$, we have \beqn
 [X_1,\ldots , [X_k, Z]\dots]\in \g^{\sst 0}\,. \hspace{3mm}
\eeq
We use (\ref{another_extension}) to show that $ \psi \left(
[X_1,\ldots,  [X_l, Z]\dots] \right)=\psi' \left( [X_1,\ldots,
[X_k, Z] \right)$. On the other hand $\psi (X_i)=\psi'(X_i)$ and
$\psi$ is an invertible map, therefore
\beqn
 [X_1, \ldots, [X_k, \psi (Z)-\psi' (Z)]\dots]=0\,.
\eeq But it is easy to see in local coordinates that
$[\g^{-1},Z]=0$ for $Z\in\g^k$, $k>0$, implies $Z=0$, so we get
inductively $\psi (Z)=\psi' (Z)$. $\square$

\section{Lie algebras of supermanifolds}
\subsection{General case}
Let $\cM$ and $\cN$ be smooth supermanifolds, $\cA$ and $\cB$ the algebras of functions on $\cM$ and $\cN$,
respectively, and $\g$ and $\h$ the corresponding super Lie algebras of vector fields.
 Suppose that the supermanifolds are supplied with a compatible $\Z -$grading as in
(\ref{isomorphism_functions}), which means that there exist vector bundles $V\to M$ and $W\to N$, such that
$\cM\simeq \Pi V$ and $\cN\simeq \Pi W$ in the category of $\Z -$graded manifolds. Let us also assume that
$\dim M$ and $\dim N$ are non-zero and the rank of the corresponding bundles
is greater than $1$ or $\rk V$ and $\rk W$ are both greater than $2$.

\begin{theorem}\label{Isomorphism_Z2} For any isomorphism
$\psi\colon\g \to \h$ of the super Lie algebras of vector fields on $\cM$ and $\cN$ respectively there exists
a diffeomorphism $\phi :\cM\to \cN$ such that $\psi (X)=(\phi^*)^{-1}\circ X \circ\phi^*$, where $\phi^*$ is
the isomorphism of the corresponding superalgebras of smooth functions functions induced by $\phi$.
\end{theorem}

\noindent\proof Taking into account that the isomorphism of super vector fields $\psi$ preserves the canonical filtration in
$\g$ and $\h$, determined by the correspondent maximal ad-nilpotent ideals, we obtain a unique bundle map
$V\to W$ as in Theorem \ref{Isomorphism_Z}. This bundle map induces an isomorphism $\chi$ of super Lie
algebras $\g\to\h$, which is also an isomorphism of the corresponding $\Z -$graded Lie algebras, associated to
the filtration, such that $\psi^{-1}\chi \colon \g\to \g$ has a trivial coset.

\vskip 2mm\noindent Now we use the result of Corollary
\ref{morphism_groups}.
 It is sufficient to prove that the kernel of
$\mathrm{p}\colon Aut_{\Z_2} (\g)\to Aut_{\Z} (\gr)$ consists of
automorphisms induced by super diffeomorphisms. Assume that $\psi$
belongs to the kernel of $\mathrm{p}$, that is, for each $X\in
\g^{\sst k}$,
 \beq
  \psi (X)=X +\psi_{2}(X)+\psi_4 (X)+\ldots \,, \hspace{3mm}
  \mathrm{where}\hspace{1mm} \psi_{2m}(X)\in \g^{\sst k+2m}\,.
 \eeq
Let us denote $\psi_2 (\e)$ by $Y$, then for each $X\in\g^{\sst
k}$ the following identity holds,
 \beq
  [\psi(\e), \psi (X)]=[\e +Y+\ldots, X+
  \psi_2 (X)+\ldots]= kX +[\e, \psi_2 (X)]+[Y, X]+\ldots\,,
 \eeq
where "..." are the higher degree terms. On the other hand,
 \beqn \psi[\e, X]=k\psi (X)= kX+k \psi_2 (X)+\ldots\eeq and $[\e, \psi_2
(X)]= (k+2)\psi_2 (X)$, therefore $\psi_2 (X)=-\frac{1}{2}[Y, X]$.
We exponentiate $\frac{1}{2}Y$ to a super automorphism of $\cA$ by
use of the exponential series, which obviously converges because
of the nilpotency of $Y$. Then the new automorphism
$\psi^{(1)}\colon =\psi\circ\mathrm{Ad}(\exp\frac{1}{2}Y)$ is
equal to the identity up to the 2d order, i.e. for each $X$ of the
degree $k$,
 \beqn\psi^{(1)}(X)-X\in \g^{\sst (k+4)}\,.
 \eeq
 Let us repeat
this procedure by induction (the number of steps will be certainly
finite because the Lie algebra $\gr$ is finitely-graded). Finally
we obtain a decomposition of $\psi$ into a (finite) product of
$\mathrm{Ad}(\exp Y_j)$ for vector fields $Y_j$ of degree $2j$,
which implies that $\psi$ is induced by a the pullback of super
diffeomorpism of $\cM$. $\square$

\subsection{Exceptional low rank cases}
Now let us consider the exceptional case, when $\rk V=1$ or $\dim
M=0$ and $\rk V\le 2$. Whatever isomorphism of the form
(\ref{isomorphism_functions}) is chosen,  $\g_0=\g^0$,
$\g_1=\g^{-1}\oplus\g^1$ such that $\g^{\pm 1}$ can be identified
with the spaces of smooth sections of certain vector bundles over
$M$. In particular, $\g^{-1}$ is always isomorphic to $\Gamma (V)$
and
 \beq\label{g1}
  \g^1 \simeq\left\{
  \be{ccc}
   \Gamma (V^*\otimes TM)\,, &\mathrm{if}\,& \rk V=1\\
   \Gamma (\Lambda^2 V^*\otimes V)\,, &\mathrm{if} \,& \dim M=0\,, \,\rk V=2
  \ee
  \right.
 \eeq

\noindent Only in this situation the canonical ideal $\g'$ is
zero, thus if $\g$ is isomorphic to the super Lie algebra $\h$ of
vector fields on another supermanifold $\cN$, the manifold $\cN$
has to satisfy the same conditions if the ranks are concerned.

\begin{lemma}\label{pm1}
 The vector fields $\pm\e$ are the only elements of the form
$f\e$, where $f$ is a smooth function on $M$, the restriction of
the adjoint action of which to $\g_1$ has only eigenvalues $\pm
1$.
\end{lemma}
\proof Indeed, for each $Y\in\g^{\pm 1}$ one has $[f\e, Y]=\pm Y$.
It is obviously true for $\g^{-1}$, whatever the base manifold $M$
is taken, and for $\g^1$ in the case of $\dim M=0$. Let us
consider the remaining case of $\g^1$ when $\dim M>0$ and $\rk
V=1$. Suppose we are given a local coordinate chart with
coordinates $\{x^i, \xi\}$ where the only $\xi$ is odd. Then the
restriction of $f\e$ and $Y$ to the local chart is
$f(x)\xi\pr_\xi$ and $\xi Z$, respectively, where $Z$ is a local
vector field on $M$. Now the simple computation gives
\beqn
 [f\e, Y]=[f\e,\xi Z]=f\xi Z+\xi^2 [f\e, Z]=fY\,,
\eeq
which finishes the proof of lemma. $\square$

\begin{proposition}\label{Except} Let $\psi\colon \g \to\h$
be a $\Z_2-$graded
 isomorphism of the Lie algebras of vector fields on $\Pi V$ and $\Pi
W$, respectively. Then either $\psi$ is induced by an isomorphism
of vector bundles $\phi \colon V\to W$ (as in Theorem
\ref{Isomorphism_Z}) or $\psi$ is the composition of two
isomorphisms, the one of which is as above and the second one is
uniquely determined by a bundle isomorphism $V\to V^*\otimes TM$
($\rk V=1$) or $V\to\Lambda^2 V^*\otimes V$ ($\dim M=0$ and $\rk
V=2$).
\end{proposition}
\proof Apparently, $\psi (\g^0)=\g^0$.
Following the similar ideas as in Theorem \ref{Isomorphism_Z}, we
immediately prove that $\psi (f\e)=\lambda (h^*)^{-1}(f)\e$
for each smooth function $f$ on $M$ where $\lambda$ is some non-zero
constant and $h\colon M\to N$
is a diffeomorphism. Using Lemma \ref{pm1}, we get $\la=\pm 1$.

\vskip 2mm\noindent If $\la=1$ then
we are in the situation of Theorem \ref{Isomorphism_Z}, which means that
$\psi$ is implemented by a bundle isomorphism $V\to W$.

\vskip 2mm\noindent Suppose $\la =-1$, then $\psi$ exchanges
$\g^{\pm 1}$ and $\h^{\mp 1}$. Applying the same argument as in
Theorem \ref{Isomorphism_Z}, we conclude that the restriction of
$\psi$ to $\g^{-1}$ is induced by a bundle map covering the
diffeomorphism $h$. In particular, if $\rk V=1$ then the
restriction of $\psi$ to $\g^{-1}$ is induced by a bundle
isomorphism $V\to W^*\otimes TN$, which implements the dimension
property $\dim N=\dim M=1$ because of the rank argument. If $\rk
V=2$ and $\dim M=0$, we obtain (in the same way) a non-degenerate
linear map $V\to \Lambda^2 W^*\otimes W$.

\vskip 2mm\noindent In both acceptable cases, when $\dim M$ equals
to $1$ or $0$, $\g^1$ and $\g^{-1}$ are isomorphic as vector
bundles on $M$. Combining any bundle isomorphism $\g^{-1}\to\g^1$,
which covers the identity diffeomorphism $M\to M$, with $\psi$, we
get a bundle isomorphism $\phi : V\to W$ which covers $h$. In
particular, this implies that $\cM$ and $\cN$ are diffeomorphic as
smooth supermanifolds. Now we can decompose $\psi$ as
$\phi_*\circ\psi_0$ where $\phi_*$ is the isomorphism of vector
fields induced by the diffeomorphism $\phi$ and $\psi_0$ is an
automorphism of $\g$ which replaces $\e$ with $-\e$, thus $\g^1$
with $\g^{-1}$.

\vskip 2mm\noindent As we have seen above, $\psi_0$ inspires a bundle map
$V\to V^*\otimes TM$ if $\dim M=1$ and $V\to \Lambda^2 V^*\otimes V$ if
$\dim M=0$. Since the representation of $\g_0$ in sections of $V$ is faithful and
any bundle isomorphism $V\to V$, commuting with the adjoint action of
all sections of $\g_0$, is the identity, we conclude that $\psi_0$ is uniquely fixed
by its restriction to $\g^{-1}$ or by a non-degenerate section of either
$(V^*)^{\otimes 2}\otimes TM$ or $\Lambda^2 V^*\otimes V^*\otimes V$ (depending on
the dimension of $M$).
$\square$

\medskip\noindent
A general corollary independent on the rank of the bundles in
question is now the following.

\begin{cor}
Two supermanifolds are diffeomorphic if and only if their super
Lie algebras of vector fields are isomorphic.
\end{cor}

\vskip1cm
\noindent Janusz GRABOWSKI\\ Polish Academy of Sciences, Institute of
Mathematics\\ \'Sniadeckich 8, P.O. Box 21\\ 00-956 Warsaw,
Poland\\Email: jagrab@impan.pl \\

\noindent Alexei KOTOV\\ University of Luxembourg, Mathematics Laboratory\\
avenue de la Fa\"{\i}encerie, 162A \\
L-1511 Luxembourg City, Grand-Duchy of Luxembourg\\Email: alexei.kotov@uni.lu\\

\noindent Norbert PONCIN\\University of Luxembourg,Mathematics
Laboratory\\avenue de la Fa\"{\i}encerie, 162A\\
L-1511 Luxembourg City, Grand-Duchy of Luxembourg\\Email: norbert.poncin@uni.lu

\end{document}